\documentclass{article}

\usepackage{amsmath, amsthm, amsfonts}
\usepackage{color}

\usepackage{hyperref}

\usepackage{longtable}

\usepackage{csquotes}

\usepackage{mdframed}
\theoremstyle{definition}
\newmdtheoremenv{boxProb}{Problem}
\newmdtheoremenv{boxDef}{Definition}
\newmdtheoremenv{boxCor}{Corollary}
\newmdtheoremenv{boxThm}{Theorem}
\newmdtheoremenv{compjob}{Computational Job}
\newmdtheoremenv{reqi}{Requirement}

\usepackage{graphicx}
\usepackage{subfigure}          %% Unterabbildungen

\usepackage{placeins}

\usepackage{tabularx}

\usepackage{lmodern,textcomp}

\usepackage{algorithm}
\usepackage{algpseudocode}

\usepackage{enumitem}

\usepackage{xspace} 
\newcommand\largeparbreak{\par\bigskip}

\newcommand{\blambda}{{\boldsymbol{\lambda}}\xspace}
\newcommand{\bxi}{\boldsymbol{\xi}\xspace}

\renewcommand{\t}{^\textsf{T}\xspace}

\newcommand{\away}[1]{}

\newcommand{\R}{\mathbb{R}\xspace}

\newcommand{\N}{\mathbb{N}\xspace}

\newcommand{\cL}{\mathcal{L}\xspace}

\newcommand{\bx}{\mathbf{x}\xspace}

\newcommand{\bei}[1]{{\mathbf{e}}\xspace}

\newcommand{\bO}{\mathbf{0}\xspace}

\newcommand{\tblambda}{\tilde{\boldsymbol{\lambda}}\xspace}

\usepackage{todonotes}

\title{Note on the Modifed Augmented Lagrangian Method for Minimization of Functions with Large Quadratic Penalties}

\author{Martin Neuenhofen}

\begin{document}

\maketitle

\begin{abstract}
In a recent work (arXiv-DOI: 1804.08072v1) we introduced the Modified Augmented Lagrangian Method (MALM) for the efficient minimization of objective functions with large quadratic penalty terms.
From MALM there results an optimality equation system that is related to that of the original objective function. But, it is numerically better behaved, as the large penalty factor is replaced by a milder factor.

In our original work, we formulated MALM with an inner iteration that applies a Quasi-Newton iteration to compute the root of a multi-variate function. In this note we show that this formulation of the scheme with a Newton iteration can be replaced conveniently by formulating a well-scaled unconstrained minimization problem.

In this note, we briefly review the Augmented Lagrangian Method (ALM) for minimizing equality-constrained problems. Then we motivate and derive the new proposed formulation of MALM for minimizing unconstrained problems with large quadratic penalties. Eventually, we discuss relations between MALM and ALM.
\end{abstract}

%\tableofcontents

\noindent
\subsubsection*{Brief summary}
We propose the iteration in Section~\ref{sec:4.1}\,. The iteration is used to solve problem \eqref{eqn:UPNP}. The iteration is also meaningful and well-behaved when $\omega_E \rightarrow +0$.

\section{Introduction}
ALM is an iterative method that is made for solving equality-constrained problems \cite[Chapter~17]{NumOpt}. The benefit is that ALM solves the original problem by minimizing a sequence of unconstrained functions that have mild quadratic penalties.

In \cite{MALM} we introduced an approach to using ALM for the solution of unconstrained functions with large quadratic penalties. This scheme we called MALM. MALM provides the same benefit as ALM. Namely, MALM solves the original problem by minimizing a sequence of unconstrained functions that have mild quadratic penalties. Hence, the large quadratic penalties are replaced by smaller ones. This is beneficial for convergence.

In our original presentation, MALM was formulated with a Quasi-Newton system: A new iterate was constructed by computing the root of a multi-variate function. In this work, through careful analysis of this multi-variate function, we are able to characterize the root of this function as a critical point for an unconstrained minimization problem. This allows us to formulate the iteration of MALM in a way that is more similar to the iteration of ALM. Hence, we are able to compare both ALM and MALM in a more elegant and meaningful way. Also, the new formulation makes the iteration of MALM simpler, which is desirable when attempting to incorporate it into other existing optimization methods.

\section{Background}
Consider the \textit{equality-constrained nonlinear program}:
\begin{equation}
	\tag{ECNP}\label{eqn:ECNP}
	\begin{aligned}
		&\operatornamewithlimits{min}_{\bx \in \R^n} 	&\quad 	f(\bx)& 	\\
		&\text{subject to} 						&		c(\bx)&=\bO \in \R^m
	\end{aligned}
\end{equation}
with the Lagrangian function $\cL(\bx,\blambda) := f(\bx)- \blambda\t \cdot c(\bx)$.

ALM is an iterative method for solving \eqref{eqn:ECNP}. The method is given a sufficiently small penalty parameter $\omega >0$, an initial guesses $\bx_0 \in \R^n$ for the primal solution, and an initial guess $\blambda_0 \in \R^m$ for the dual solution. The iteration of ALM then generates sequences $\lbrace\bx_k\rbrace_{k \in \N_0}$, $\lbrace\blambda_k\rbrace_{k\in\N_0}$, where for $k \in \N$ it holds:
\begin{enumerate}
	\item $\bx_k$ minimizes
	\begin{align*}
		f_k(\bx) := \cL(\bx,\blambda_{k-1}) + \frac{1}{2 \cdot \omega} \cdot \|c(\bx)\|_2^2
	\end{align*}
	and is computed with a numerical method for unconstrained minimization, using the initial guess $\bx_{k-1}$\,.
	\item After $\bx_k$ has been computed, the dual vector $\blambda_k$ is computed as
	\begin{align*}
		\blambda_k := \blambda_{k-1} - \frac{1}{\omega} \cdot c(\bx_k)\,.
	\end{align*}
\end{enumerate}
Under mild conditions on the problem \eqref{eqn:ECNP} and a sufficiently small choice for $\omega$, that however is finite and does not need to converge to zero, the iteration of ALM converges to a stationary point $\bx_\infty,\blambda_\infty$ of the following optimality equations due to Karush, Kuhn and Tucker (KKT) \cite{NumOpt,Lancelot}:
\begin{align*}
	\nabla_\bx \cL(\bx,\blambda) &= \bO\in\R^n\,,\\
	\nabla_\blambda \cL(\bx,\blambda) &= \bO\in\R^m\,.
\end{align*}

\section{Motivation of MALM}
Sometimes, instead of \eqref{eqn:ECNP}, we may rather want to minimize the following problem, that we call \textit{unconstrained penalty nonlinear program}:
\begin{equation}
	\tag{UPNP}\label{eqn:UPNP}
	\begin{aligned}
		&\operatornamewithlimits{min}_{\bx \in \R^n} 	&\quad 	f(\bx) + \frac{1}{2 \cdot \omega_E} \cdot \|c(\bx)\|_2^2
	\end{aligned}
\end{equation}
where typically $\omega_E>0$ is a very small user-defined parameter.

There are various reasons why, for given functions $f,c$, one may rather prefer solving \eqref{eqn:UPNP} than \eqref{eqn:ECNP}.
\begin{enumerate}
	\item Problem \eqref{eqn:UPNP} is always feasible. Thus, we do not run into any numerical issues related to infeasibility. This is beneficial when the user simply asks for a solution under any circumstances.
	\item It is impossible to verify numerically whether \eqref{eqn:ECNP} is feasible. This is so already because of the reason that in all but very few cases the equations $c(\bx)=\bO$ do not have a digitally representable solution for $\bx$.
	\item Unless assuming a constraint qualification (e.g. LICQ), the optimality equations of \eqref{eqn:ECNP} are ill-posed, and hence do not admit a meaningful numerical treatment. In contrast, \eqref{eqn:UPNP} is always well-posed, with a condition number that can be controlled by the parameter $\omega_E>0$.
	\item There are applications where $m>n$. This immediately causes violation of all constraint qualifications. However, \eqref{eqn:ECNP} may still possess feasible points in such cases. In these situations, it constitutes a sound approach to generate numerical solutions of \eqref{eqn:ECNP} by solving \eqref{eqn:UPNP} for very small values of $\omega_E>0$.
	\item In general, the numerical objective of satisfying the constraints $c$ is in contradiction to the numerical objective of minimizing $f$. Hence, from a numerical point of view, it is impossible to think of any method that does not actually minimize a bias (merit function) between $f$ and a norm of $c$. In \eqref{eqn:UPNP}, at least we know this merit function and can analyze its meaning with respect to an application.
	\item In logical consequence of the former bullet, many numerical methods for treating \eqref{eqn:ECNP} are actually based on minimizing \eqref{eqn:UPNP} for decreasing values of $\omega_E>0$. The methods described in the following literature serve as examples: \cite{Amudson,ChenGoldfarb,ForsgrenGill}.
	\item Finally and not yet mentioned, there simply exist problems where people literally wish to minimize a function of the form \eqref{eqn:UPNP}. Examples arise from applications where $n,m$ are dimensions of discretized infinite spaces for $\bx,\blambda$. This is the case, e.g., in constrained optimization with partial differential equations. It is for instance in these situations, where it can make sense to consider discretizations where $m>n$, and treat the resulting problem as \eqref{eqn:UPNP}. To give an example, in \cite{myOCP} we present a discretization for optimal control problems that results in nonlinear programs of the form \eqref{eqn:UPNP}.
\end{enumerate}

If we were to minimize \eqref{eqn:UPNP} with a numerical method for unconstrained minimization, then there would be an issue. This issue would be, that the large quadratic penalty terms in the objective impede the rate of convergence of the numerical method. Hence, it seems desirable to apply a technique such as ALM, that manages to solve a minimization problem with only mild quadratic penalties; i.e., it only uses quadratic penalties with values $\omega \gg \omega_E$. Usually, these milder choices of $\omega$ form well-scaled objectives and hence promotive fast convergence.

However, it is unclear how ALM could be applied to treat \eqref{eqn:UPNP}. This is because after all, \eqref{eqn:UPNP} and \eqref{eqn:ECNP} are two entirely different problem statements.

\largeparbreak

In \cite{MALM}, we made an approach to reformulating \eqref{eqn:UPNP} as an instance of \eqref{eqn:ECNP}, by introducing an auxiliary variable $\bxi \in \R^m$:
\begin{subequations}
	\label{eqn:Subst}
	\begin{align}
		&\operatornamewithlimits{min}_{\bx \in \R^n\,,\bxi \in \R^m} 	&\quad 	f(\bx) + \frac{\omega_E}{2} \cdot \|\bxi\|_2^2& \label{eqn:SubstObjective}	\\
		&\text{subject to} 						&		c(\bx) + \omega_E \cdot \bxi &=\bO \in \R^m \label{eqn:SubstConstraints}
	\end{align}
\end{subequations}
This is advantageous, because the new objective \eqref{eqn:SubstObjective} is well-scaled in $\bx,\bxi$ even for the smallest values of $\omega_E$, and the constraints \eqref{eqn:SubstConstraints} are well-scaled in $\bx,\bxi$ as well and feasible for all values of $\omega_E>0$.

\largeparbreak

In \cite{MALM}, after formulating \eqref{eqn:Subst}, we then applied ALM to \eqref{eqn:Subst}. We arrived at an optimality system that was linear in $\bxi$. We eliminated $\bxi$ from that system. We were then able to establish an iterative scheme, that we called MALM. This scheme is effectively the iteration of ALM applied to \eqref{eqn:Subst}.

The iteration goes as follows. Provided are initial guesses $\bx_0,\blambda_0$, and a sufficiently small value $\omega \gg \omega_E$. The iteration of MALM then generates sequences $\lbrace\bx_k\rbrace_{k \in \N_0}$, $\lbrace\blambda_k\rbrace_{k\in\N_0}$, where for $k \in \N$ it holds:
\begin{enumerate}
	\item $\bx_k$, together with an auxiliary variable $\tblambda \in \R^m$, solves the following equation system
	\begin{align*}
		F_k(\bx,\tblambda) := \begin{pmatrix}
			\nabla_\bx \cL(\bx,\blambda_{k-1} + \tblambda)\\
			c(\bx) + \omega_E \cdot \blambda_{k-1} + (\omega + \omega_E) \cdot \tblambda
		\end{pmatrix}\,.
	\end{align*}
	In our original work \cite{MALM}, we proposed that the root $(\bx_k,\tblambda)$ of $F_k$ be computed using a particular Quasi-Newton line-search method. To this end, we presented both a suitable inertia-correction scheme for the step-direction obtained from the Newton system and a suitable merit-function for which this step-direction is a descent-direction. We referred to literature results that proved global and fast local convergence of this scheme.
	\item After $\bx_k,\tblambda$ has been computed, the dual vector $\blambda_k$ is computed as
	\begin{align*}
	\blambda_k := \blambda_{k-1} + \tblambda
	\end{align*}
\end{enumerate}

\section{Associated objective for MALM}
The iteration of MALM is a superior approach over minimizing \eqref{eqn:UPNP} with an iteration for unconstrained minimization. But, MALM has one detail, that makes it un-elegant and difficult to handle. We point out this detail in the following:

The iteration of ALM defines $\bx_k$ as the local unconstrained minimizer of $f_k$. In contrast to that, the iteration of MALM defines $\bx_k$ as the root of a particular optimality system. It would be nice if instead we could define $\bx_k$ in MALM as the unconstrained minimizer of a function like $f_k$, as well.

\largeparbreak

The proposal of this note is that we found a suitable function for $f_k$. Namely, we can interpret $F_k$ as the optimality system for the following unconstrained minimization problem:
\begin{equation}
	\tag{UMP$_k$}\label{eqn:UMP}
	\begin{aligned}
		&\operatornamewithlimits{min}_{\bx \in \R^n} 	&\quad 	f_k(\bx) := \cL(\bx,\blambda_k) + \frac{1}{2 \cdot (\omega+\omega_E)} \cdot \|c(\bx)+\omega_E \cdot \blambda_{k-1}\|_2^2
	\end{aligned}
\end{equation}
This minimality problem we found as follows: We eliminate the second row in $F_k$ and obtain for $\tblambda$ the formula
$$ 	\tblambda = \frac{-1}{\omega + \omega_E} \cdot \big(\,c(\bx_k) + \omega_E \cdot \blambda_{k-1}\,\big) 	$$
in terms of $\bx_k,\blambda_{k-1}$. We insert it into the first row of $F_k$. The first row of $F_k$ with $\tblambda$ inserted is
$$ 	\nabla_\bx \cL(\bx_k,\blambda_{k-1}) - \nabla c(\bx_k) \cdot \frac{-1}{\omega + \omega_E} \cdot \big(\,c(\bx_k) + \omega_E \cdot \blambda_{k-1}\,\big) \,.	$$
It turns out that this expression coincides with the gradient of our proposed choice for $f_k$. Hence, minimizers of $f_k$ are roots of $F_k$.

\subsection{Proposed new form for MALM}\label{sec:4.1}
The formulation of $f_k$ allows us to simplify our initial proposal for MALM into the following new iterative scheme, that is described below.

The iteration goes as follows. Provided are initial guesses $\bx_0,\blambda_0$, and a sufficiently small value $\omega \gg \omega_E$. The iteration of MALM then generates sequences $\lbrace\bx_k\rbrace_{k \in \N_0}$, $\lbrace\blambda_k\rbrace_{k\in\N_0}$, where for $k \in \N$ it holds:
\begin{enumerate}
	\item $\bx_k$ minimizes
	\begin{align*}
	f_k(\bx) := \cL(\bx,\blambda_{k-1}) + \frac{1}{2 \cdot (\omega+\omega_E)} \cdot \|c(\bx)+\omega_E \cdot \blambda_{k-1}\|_2^2
	\end{align*}
	and is computed with a numerical method for unconstrained minimization, using the initial guess $\bx_{k-1}$\,.
	\item After $\bx_k$ has been computed, the dual vector $\blambda_k$ is computed as
	\begin{align*}
	\blambda_k := \blambda_{k-1} - \frac{1}{\omega+\omega_E} \cdot \big(\,c(\bx_k)+\omega_E \cdot \blambda_{k-1}\,\big)\,.
	\end{align*}
\end{enumerate}
The second step remains unchanged from the originally proposed MALM. We just inserted the above expression for $\tblambda$ in order to vanish it and clean up the iteration.

\subsection{Discussion and benefits of the result}
We find that our proposed new MALM is a direct generalization of ALM. This is because for $\omega_E=0$ the two schemes coincide.

The benefit of our method is that it is both capable of solving \eqref{eqn:ECNP} and \eqref{eqn:UPNP}, while ALM can only treat the latter. If we choose $\omega_E>0$, then MALM solves \eqref{eqn:UPNP} for that respective value of $\omega_E$. If instead we choose $\omega_E=0$, then MALM coincides with ALM and solves \eqref{eqn:ECNP}. Hence, for every $\omega_E \in [0,\infty)$ the iteration is well-defined. Therefore, MALM offers a convenient treatment for solving problems of both classes, and for solving problems \eqref{eqn:UPNP} with literally arbitrarily small values of $\omega_E$.

From a numerical perspective, the proposed new MALM iteration appears as suitable for computations as does ALM. We justify this: In ALM, for the update of the Augmented Lagrangian in step 2, there is a division through a small parameter $\omega$. In the new proposed MALM, this divisor is augmented with $\omega_E$. Certainly, this augmentation does not introduce any numerical difficulty, because it only increase of the divisor. Hence, we believe that MALM resembles the desirable practical convergence behavior of ALM, while being applicable to a more general setting.

\section{Summary}
We have modified the Augmented Lagrangian Method in a small detail. Namely, we added a parameter $\omega_E \in [0,\infty)$ into the iteration formulas. This modification allows us to treat problems of either form \eqref{eqn:ECNP} or \eqref{eqn:UPNP}.

Our Modified Augmented Lagrangian Method is numerically as well-behaved as the original Augmented Lagrangian Method. This relates to numerical well-posedness and to the rate of convergence. Our modified method is insensitive to small choices of $\omega_E>0$ and coincides with the original Augmented Lagrangian in the limit $\omega_E = 0$.

\largeparbreak

Further research may be related to the following subject. We understand that --\,whenever existent\,-- a solution of \eqref{eqn:ECNP} can be approached by solving a sequence of \eqref{eqn:UPNP} for $\omega_E \searrow +0$\,. Now, with our Modified Augmented Lagrangian Method, we even have a convenient tool at hand that can solve \eqref{eqn:UPNP} in the limit $\omega_E = 0$\,. This motivates research investigations in the following direction: What solution does the sequence \eqref{eqn:UPNP} for $\omega_E \searrow +0$ converge to when \eqref{eqn:ECNP} does not have any solutions. And, related to this situation, to which values does the (Modified) Augmented Lagrangian method with $\omega_E \searrow +0$ converge to? Or, does the limit exist in these cases?

Certainly, a progress in this regard is provided in our analysis of a method for convex quadratic programming. We show that for $f$ convex quadratic and $c$ linear the solution of \eqref{eqn:UPNP} for $\omega_E \searrow +0$ converges to a solution of the bi-objective minimization problem
\begin{equation*}
	\begin{aligned}
		&\operatornamewithlimits{min}_{\bx \in \R^n} 	&\quad 	f(\bx)&\,,\\
		&\text{subject to} 								&\quad \|c(\bx)\|_2& \text{ is locally minimal\,.}
	\end{aligned}
\end{equation*}
We conjecture that a similar result can also be verified for the more general case of non-convex nonlinear programming.

\largeparbreak

Finally, justified by our proposed method, it would be desirable if people could please provide optimization software that is capable of treating problems \eqref{eqn:UPNP}, and that does not return on error when $m>n$\,.

\FloatBarrier

%\nocite{*}
\bibliography{NOTE_MALM_bib}

\begin{thebibliography}{1}

\bibitem{Amudson}
P.~Armand, J.~Benoist, R.~Omheni, and V.~Pateloup.
\newblock Study of a primal-dual algorithm for equality constrained
  minimization.
\newblock {\em Comput. Optim. Appl.}, 59(3):405--433, 2014.

\bibitem{ChenGoldfarb}
L.~Chen and D.~Goldfarb.
\newblock Interior-point l2-penalty methods for nonlinear programming with
  strong global convergence properties.
\newblock {\em Mathematical Programming}, 108(1):1--36, Aug 2006.

\bibitem{Lancelot}
A.~R. Conn, N.~I.~M. Gould, and P.~L. Toint.
\newblock {\em Lancelot: A Fortran Package for Large-Scale Nonlinear
  Optimization (Release A)}.
\newblock Springer Publishing Company, Incorporated, 1st edition, 2010.

\bibitem{ForsgrenGill}
A.~Forsgren and P.~E. Gill.
\newblock Primal-dual interior methods for nonconvex nonlinear programming.
\newblock {\em SIAM J. Optim.}, 8(4):1132--1152, 1998.

\bibitem{myOCP}
M.~P. {Neuenhofen}.
\newblock High-order convergent finite-elements direct transcription method for
  constrained optimal control problems.
\newblock {\em ArXiv e-prints. DOI:1712.07761}, dec 2017.

\bibitem{MALM}
M.~P. {Neuenhofen}.
\newblock Modified augmented lagrangian method for the minimization of
  functions with quadratic penalty terms.
\newblock {\em ArXiv e-prints. DOI:1804.08072}, apr 2018.

\bibitem{NumOpt}
J.~Nocedal and S.~J. Wright.
\newblock {\em Numerical optimization}.
\newblock Springer Series in Operations Research and Financial Engineering.
  Springer, New York, second edition, 2006.

\end{thebibliography}
\bibliographystyle{plain}

\end{document}